\documentclass[10pt]{article}
\usepackage{latexsym, amsfonts,amsthm,amsmath,amssymb,amscd,verbatim,multicol,tabularx,graphicx}

\newtheorem{theorem}{Theorem}[section]

\newtheorem{proposition}[theorem]{Proposition}

\theoremstyle{definition}
\newtheorem{definition}[theorem]{Definition}
\newtheorem{remark}[theorem]{Remark}

\newcommand{\Ker}{\mathop{\rm Ker}\nolimits} % kernel
\newcommand{\Ima}{\mathop{\rm Im}\nolimits}   % image

\def\G{\Gamma}
\def\l{\lambda}
\def\A{\mathcal A}
\def\Hr{\widetilde{H}}
\def\Cr{\widetilde{C}}
\def\pr{\widetilde{p}}
\def\Rr{\widetilde{R}}
\newcommand{\R}{\mathbb R} % real numbers
\newcommand{\Z}{\mathbb Z} % real numbers
\newcommand{\M}{\mathbb M} % module over \A
\newcommand{\ig}{\includegraphics}

\begin{document}

\title{Reduced chromatic graph cohomology}
\author{Michael Chmutov, Elena Udovina}
%\date{}
\maketitle

\begin{abstract}
In this paper we give a new characterization of the $h$-vector of the 
chromatic polynomial of a graph, i.e. the vector $(h_0, \dots, h_n)$ 
of coefficients of the chromatic polynomial 
$$p_\G(\l) = h_0\l(\l-1)^{n-1} - h_1\l(\l-1)^{n-2} + \dots + (-1)^{n-1}h_{n-1}\l.$$
We introduce reduced chromatic cohomology of a graph and show that 
$h_i$ are its Betti numbers. We then discuss various combinatorial 
properties of these cohomologies.
\end{abstract}

%---------------
%---------------
%0. Introduction
%---------------
%---------------

\section*{Introduction} \label{s:intro}
In \cite{HGR} L.~Helme-Guizon and Y.~Rong introduced a bigraded cohomology 
theory for graphs whose graded Euler characteristic is equal to the chromatic 
polynomial. Their work was motivated by the development of the Khovanov 
cohomology in knot theory \cite{Kho}. The theory also suggests a 
notion of reduced Khovanov cohomology. The reduced cohomology for graphs was 
introduced by J.Przytycki \cite{Pr}.

In this paper we work with a specialization of the reduced cohomology to the 
algebra $\A=\R[x]/(x^2)$ and the module $\M$ over $\A$ being the ideal of $\A$ 
generated by $x$ (we shift the degree by $-1$ as compared to \cite{Pr} for 
better agreement with combinatorial formulas for the chromatic polynomial). 
Our results include the following. 
The reduced chromatic cohomology groups are concentrated on one diagonal 
(Proposition \ref{diagonal}). 
If an edge $e$ is not a bridge then $(i,j)$-th reduced cohomology of $\G$ 
is a direct sum of $(i,j)$-th reduced cohomologies of $\G/e$ and $\G-e$ 
(Proposition \ref{direct sum}).
Let $p_\G(\l)$ denote the chromatic polynomial of $\G$. Then the graded Euler 
characteristic of reduced cohomologies $\sum_{i,j} (-1)^i q^j \dim(\Hr^{i,j}(\G))$
is equal to $p_\G(1+q)/(1+q)$ (Proposition \ref{euler chromatic}).
The reduced cohomology of a one vertex union of two graphs is equal to the tensor product 
of the reduced cohomologies of its factors (Proposition \ref{vertex union}). 
In section \ref{s:polynomial} we describe the relation between the standard cohomology (over a field) and 
the reduced cohomology. In particular, it implies that the reduced cohomologies are determined 
by the chromatic polynomial. As a consequence of this we can conclude that the reduced 
cohomologies depend only on the matroid type of the graph. However we prefer to give a separate 
proof (in Section \ref{s:matroid}) of this fact using Whitney twists. 

This work was done during the Summer'05 VIGRE working group ``Knot theory and Combinatorics" 
at the Ohio State University funded by NSF grant DMS-0135308. We are grateful to the participants 
of the working group and Prof.~J.Przytycki for the useful discussions.
%----------------------
%----------------------
%0. Preliminary Results
%----------------------
%----------------------

\section{Definitions and preliminary results} \label{s:def}

\begin{figure}[bt]
\label{f:prelim}
\centerline{\ig[width=340pt] {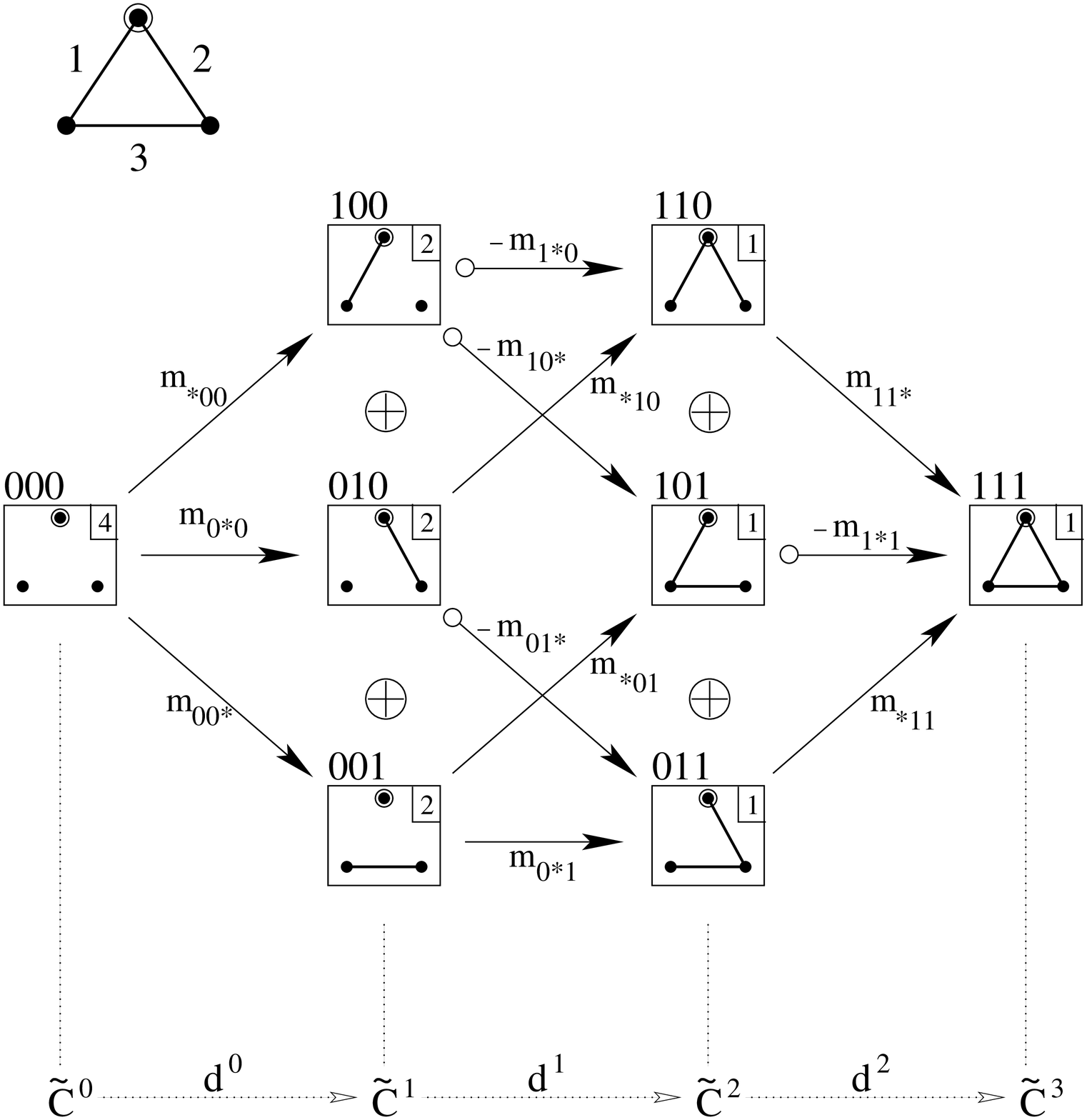}}
\caption{A ``smoothing'' diagram for a triangle.}
\end{figure}

For a graph $\G$, a {\it state} is a spanning subgraph of $\G$, that is a 
subgraph of $\G$ containing all the vertices of $\G$ and a subset of the 
edges. The number of edges in a state is called its {\it dimension}. Choose a 
vertex of $s$ as a base point. A {\it reduced enhanced state} $S$ is a state 
whose connected components are colored in two colors, $x$ and $1$, 
and the component with the base point is always colored in $x$. The {\it degree} of $S$
is the number of connected components 
colored in $x$ minus $1$. 
The {\it cochain group} 
$\Cr^{i,j}$ is defined to be the real vector space spanned by all enhanced 
states of dimension $i$ and degree $j$. These notions are illustrated in 
Figure \ref{f:prelim}, with the base point circled. This picture is 
similar to Bar-Natan's \cite{BN}.

Here every square box represents a vector space spanned by all reduced enhanced 
states with the indicated underlying state, and its dimension is shown in the upper 
right hand corner. The direct sum of these vector 
spaces located in the $i$-th column gives the cochain group 
$\Cr^i=\bigoplus\limits_j \Cr^{i,j}$.
The boxes are labeled by strings of $0$'s and $1$'s which encode the 
edges participating in the corresponding states. To turn the cochain groups 
into a {\it cochain complex} we define a differential 
$d^{i,j}:\Cr^{i,j}\to \Cr^{i+1,j}$.
On a vector space corresponding to a given state (box) the differential can 
be defined as adding an edge to the corresponding state in all possible 
ways, and then coloring the connected components of the obtained state 
according to the following rule. Suppose we are adding an edge $e$ to a reduced 
enhanced state $S$. Then, if the number of connected components is not 
changed, we preserve the same coloring of connected components of the new 
state $S\cup e$. If $e$ connects two different connected components of $S$, 
then the color of the new component of $S\cup e$ is defined by the multiplication
$$1\times 1:=1,\qquad 1\times x:=x,\qquad x\times 1:=x,\qquad x\times x:=0\ .$$
In some cases, if the number of edges of $S$ whose index is less than that of $e$ is odd, 
we should take the target reduced enhanced state $S\cup e$ with the coefficient $-1$. 
These are shown in the picture above by arrows with little circles at their tails. 

Since $d$ is defined in the same way as for standard cohomologies, we have
$d^{i+1,j}\circ d^{i,j}=0$ converting our cochain groups into a bigraded cochain 
complex $\Cr^{*,*}(\G)$. We call its cohomology groups
the {\it reduced chromatic cohomology} of the graph $\G$:
$$\Hr^{i,j}(\G)\ :=\ 
\frac{\Ker(d: \Cr^{i,j}(\G)\to \Cr^{i+1,j}(\G))}{
      \Ima(d: \Cr^{i-1,j}(\G)\to \Cr^{i,j}(\G))}\ .
$$

\begin{remark} Reduced cohomologies are independent of the ordering on 
edges since the isomorphism between cochain groups coming from different 
edge orderings in \cite[Theorem 14]{HGR} works verbatim in this case.
\end{remark}

\begin{remark}\label{components} The cohomology groups of the graph are 
tensor products of the cohomology groups of the connected components. 
Therefore, in all that follows we consider only connected graphs.
\end{remark}

\begin{remark} The long exact sequence of cohomology groups can be adapted 
for the reduced cohomologies. If $e$ is an edge of $\G$ and $e$ is not a bridge 
then the differential $d$ commutes with the 
maps from the short exact sequence of complexes
$$
0 \rightarrow \Cr^{i-1}(\G/e) \rightarrow \Cr^i(\G) \rightarrow \Cr^i(\G - e) \rightarrow 0,
$$
in the same way as for the non-reduced cochain complexes and hence the long 
exact sequence is analogous:
$$
0 \rightarrow \Hr^0(\G) \rightarrow \Hr^0(\G - e) \rightarrow \Hr^{0}(\G/e) 
  \rightarrow \Hr^1(\G) \rightarrow\dots.
$$
However, if removal of $e$ separates $\G$ into two components 
$\G_1$ and $\G_2$ with $\G_1$ containing the chosen vertex, then the short exact 
sequence of complexes is 
$$
0 \rightarrow \Cr^{i-1}(\G/e) \rightarrow \Cr^i(\G) \rightarrow \Cr^i(\G_1)\otimes C^i(\G_2) \rightarrow 0.
$$
The corresponding long exact sequence in cohomologies is 
$$
0 \rightarrow \Hr^0(\G) \rightarrow \Hr^0(\G_1)\otimes H^0(\G_2) \rightarrow \Hr^{0}(\G/e) 
  \rightarrow \Hr^1(\G) \rightarrow\dots.
$$
\end{remark}

\begin{remark}\label{loop} As in \cite[Propositions 19, 20]{HGR}, 
the cohomology groups of a graph with a loop are trivial, and the 
cohomology groups of a graph with multiple edges are unchanged if 
the multiple edges are replaced by single edges. Hence, in all 
that follows, the graphs will be simple.
\end{remark}

\begin{remark}\label{tree} If $e$ is a pendant edge of 
$\G$ then 
$\Hr^{i,j}(\G) \cong \Hr^{i,j-1}(\G/e)$. 
The proof is exactly the same as in \cite[Theorem 24]{HGR} 
(notice that we use the second exact sequence for reduced cohomologies here) 
since it only used the long exact sequence of cohomologies and the fact 
that $1$ is an identity in the algebra. Note that for a single 
vertex the reduced cohomology group is simply $\R$, so for a 
tree on $n$ vertices the reduced cohomology group is $\R(q^{n-1})$, i.e. 
it is one-dimensional and concentrated in cohomological dimension $0$.
\end{remark}

The vector space corresponding to a single regular vertex without any edges 
is isomorphic to the algebra of truncated polynomials 
$\A:=\R[x]/(x^2)$, while the vector space corresponding to the chosen vertex is 
the ideal of $\A$ generated by $x$. We can generalize the construction to an arbitrary 
algebra $\A$ and an $\A$-module $\M$ as follows. We think about a box space 
of an arbitrary graph $\G$ as a tensor product of the module $\M$ with a tensor 
power of the algebra $\A$ whose tensor factors are in one-to-one correspondence
with the connected components that do not contain the chosen vertex. Then our 
multiplication rule for the differential turns out to be the multiplication operation 
$\A\otimes \A\to \A$
in the algebra $\A$ together with the multiplication in the module $\M$. This approach 
allows to generalize the definition of reduced chromatic cohomology to an arbitrary algebra 
$\A$ (see \cite{Pr} for a discussion of this approach).

%-----------------------------------
%-----------------------------------
%2. Properties of Reduced Cohomologies
%-----------------------------------
%-----------------------------------

\section{Properties of Reduced Cohomologies}\label{s:properties}

\begin{proposition}\label{diagonal}
$\Hr^{i,j}(\G)=0$ unless $i+j=n-1$, where $n$ is 
the number of vertices of $\G$.
\end{proposition}
\begin{proof}
The proof is by induction on the number of edges.

{\it Base Case.} There is only one graph with 0 edges: the one-vertex tree. The 
cohomology of a single vertex is $\R$.

{\it Induction Step.} If $\G$ is a tree, the assertion of the proposition 
follows from Remark \ref{tree}. Otherwise, let $e$ be an esde that is not a 
bridge. The relevant portion of the long exact sequence is as follows:
\[
\dotsc \rightarrow \Hr^{i-1,j}(\G/e) \rightarrow \Hr^{i,j}(\G) \rightarrow \Hr^{i,j}(\G-e) \rightarrow \dotsc
\]
Since $\G/e$ and $\G-e$ have fewer edges than $\G$, by the induction 
hypothesis $\Hr^{i-1,j}(\G/e)=0$ unless 
$i-1+j=(n-1)-1$ and $\Hr^{i,j}(\G-e)=0$ unless 
$i+j=n-1$. From exactness, $H^{i,j}(\G)=0$ unless $i+j=n-1$.
\end{proof}

\begin{proposition}\label{direct sum}
Let $\G$ be a simple connected graph with $n$ vertices. Let $e$ be an edge that is 
not a bridge. Then
$\Hr^{i,j}(\G) \cong \Hr^{i-1,j}(\G/e) \oplus \Hr^{i,j}(\G-e)$.
\end{proposition}

\begin{proof}
Note that unless $i+j=n-1$, by Proposition \ref{diagonal} all the cohomologies are 
zero. When $i+j=n-1$, the relevant segment of the long exact sequence looks as follows:
\[
0 = \Hr^{i-1,j}(\G) \rightarrow \Hr^{i-1,j}(\G/e) \rightarrow \Hr^{i,j}(\G) \rightarrow \Hr^{i,j}(\G-e) \rightarrow \Hr^{i,j}(\G/e) = 0
\]
By exactness, $\Hr^{i,j}(\G) \cong \Hr^{i-1,j}(\G/e) \oplus \Hr^{i,j}(\G-e)$.
\end{proof}

\begin{proposition}\label{independence}
Reduced cohomologies $\Hr^{i,j}(\G)$ are independent of 
the choice of the special vertex.
\end{proposition}

\begin{proof}
The proof is by induction on the number of edges.

{\it Base Case.} If $\G$ has no edges, then it is the one-vertex graph, and 
there is nothing to prove.

{\it Induction Step.} For $\G$ a tree, the proposition follows from 
remark \ref{tree}, since the cohomology groups are the same regardless of 
vertex choice. Else, let $\G$ and $\G'$ correspond to the same 
graph but with different special vertices, $v$ and $v'$ respectively. 
Let $e$ be an edge of $\G$ (also of $\G'$) which is not a bridge. 
By Proposition \ref{direct sum},
\begin{align*}
\Hr^{i}(\G ) &\cong \Hr^{i-1}(\G  / e) \oplus \Hr^i(\G  - e)\\
\Hr^{i}(\G') &\cong \Hr^{i-1}(\G' / e) \oplus \Hr^i(\G' - e)
\end{align*}
By inductive hypothesis, 
$\Hr^{i-1}(\G/e) \cong \Hr^{i-1}(\G'/e)$ 
and $\Hr^i(\G-e) \cong \Hr^i(\G'-e)$, 
since these represent the same graph but with different special 
vertices, the images of $v$ and $v'$ respectively. Hence, 
$\Hr^{i}(\G) \cong \Hr^i(\G')$.
\end{proof}
%----------------------
%----------------------
%3. Cohomologies of Union
%----------------------
%----------------------

\section{Cohomologies of Union}\label{s:union}

The motivation for introducing reduced cohomologies is the following 
property of the chromatic polynomial: if $\G$ is obtained from 
$\G_1$ and $\G_2$ by taking a vertex $v_1 \in \G_1$ and a 
vertex $v_2 \in \G_2$ and glueing them together 
($\G = \G_1 * \G_2$), then
\[
p_\G(\l) = \frac{p_{\G_1}(\l) p_{\G_2}(\l)}{\l}.
\]
Introducing the reduced polynomial $\pr_\G(\l) = \l^{-1} p_\G(\l)$, we get
\[
\pr_\G(\l) = \pr_{\G_1}(\l) \pr_{\G_2}(\l).
\]

We now establish that the reduced cohomologies form the categorification 
of the reduced chromatic polynomial and have this multiplication property.

\begin{proposition}\label{euler chromatic}
The graded Euler characteristic of the reduced cochain complex $\Cr(\G)$ is equal 
to the reduced chromatic polynomial $\pr_\G(\l)$ with $\l = 1+q$.
\end{proposition}
\begin{proof}
If $\G$ is a tree on $n$ vertices and the reduced chromatic 
polynomial is 
$\l^{-1}\bigl(\l(\l-1)^{n-1}\bigr) = (\l-1)^{n-1}=q^{n-1}$ 
while the zeroth reduced cohomology group is $\R(q^{n-1})$. By 
Proposition \ref{direct sum}, the graded Euler characteristic 
satisfies the contraction-deletion relation of the chromatic 
polynomial (which is also satisfied by the reduced chromatic 
polynomial). Induction on the number of edges completes the 
proof.
\end{proof}

\begin{proposition}\label{vertex union}
If 
$\G = \G_1 * \G_2$, then $\Hr(\G) = \Hr(\G_1) \otimes \Hr(\G_2)$.
\end{proposition}

\begin{proof}
Since the cohomology groups are independent of the choice of special 
vertex, we may suppose $\G_1$ and $\G_2$ are joined by 
identifying their special vertices; the resulting 
vertex is special in the union.

The proof is by induction on the number of edges of $\G_2$.

{\it Base Case.} If $\G_2$ has no edges, then it is the single-vertex 
graph, so $\Hr(\G) = \Hr(\G_1)$ 
(since the graphs are the same), and $\Hr(\G_2) = \R$.

{\it Inductive Step.} If $\G_2$ is a tree with $n$ vertices, 
$\Hr(\G)$ is obtained from $\Hr(\G_1)$ 
via Proposition \ref{tree} as 
$\Hr(\G) = \Hr(\G_1) \otimes \R(q^{n-1})$. 
By the same proposition, $\R(q^{n-1})$ is the cohomology 
of $\G_2$. If $\G_2$ is not a tree, let $e$ be an edge of 
$\G_2$ that is not a bridge. By 
Proposition \ref{direct sum},
\begin{align*}
\Hr(\G_2) &\cong \Hr(\G_2 / e) \oplus \Hr(\G_2 - e),\\
\Hr(\G) &\cong \Hr(\G / e) \oplus \Hr(\G - e).
\end{align*}
By inductive hypothesis,
\begin{align*}
\Hr(\G / e) &\cong \Hr(\G_1) \otimes \Hr(\G_2 / e),\\
\Hr(\G - e) &\cong \Hr(\G_1) \otimes \Hr(\G_2 - e).
\end{align*}
Taking the direct sum,
\[
\Hr(\G) \cong \Hr(\G_1) \otimes \bigl( \Hr(\G_2 / e) \oplus \Hr(\G_2 - e) \bigr) \cong \Hr(\G_1) \otimes \Hr(\G_2).
\]
\end{proof}
%---------------
%---------------
%4. Matroid Type
%---------------
%---------------

\section{Matroid Type}\label{s:matroid}

\noindent A {\it Whitney twist} on a graph $\G$ can be defined as 
follows \cite{Wh, Hug}. Let $\G_1$ and $\G_2$ be two graphs. 
Pick edges $e_1 \in \G_1$ and $e_2 \in \G_2$. Construct a new 
graph by gluing the edges $e_1 \in \G_1$ and $ e_2 \in \G_2$ 
together (with their endpoints) and then removing the resulting single 
edge from the graph. In general this can be done in two ways depending 
on how we glue $e_1$ with $e_2$. If one of them is $\G$ then the 
other is a Whitney twist of $\G$. Whitney proved that two 2-connected 
graphs have the same matroid type iff one can be obtained from the other 
by a sequence of Whitney twists.
 
We show that the reduced cohomology sequence of a graph is invariant 
under the Whitney twist. From this we derive that the reduced cohomology 
sequence is an invariant of the matroid type of the graph.

\begin{proposition}\label{twist}
If $\G$ and $\G'$ are related by a Whitney twist, 
$\Hr(\G) \cong \Hr(\G')$.
\end{proposition}

\begin{proof}
Let $G$ and $G'$ be obtained by joining $\G_1$ and $\G_2$ along 
$e$. Set $\G = G - e$; then $\G' = G' - e$ is its Whitney twist. 
By Proposition \ref{direct sum},
\[
\Hr^i(G) \cong \Hr^{i-1}(G / e) \oplus \Hr^i(G - e)
\]
and similarly for $G'$. We are interested in proving 
$\Hr^i(\G) \cong \Hr^i(\G')$, 
but it suffices to prove the isomorphisms $\Hr^i(G) \cong \Hr^i(G')$, 
and $\Hr^i(G/e) \cong \Hr^i(G'/e)$.

Note that $G / e = (\G_1 / e)*(\G_2 / e) = G' / e$; 
hence, the cohomologies in both cases are just the tensor product 
of the cohomologies of $\G_1 / e$ and $\G_2 / e$. To show 
the isomorphism of the cohomology groups of $G$ and $G'$, we induct 
on the number of edges of $\G_2$.

{\it Base Case.} $\G_2$ cannot have less than one edge, since 
we have to glue $\G_1$ and $\G_2$ together along an edge. 
If $\G_2$ has exactly one edge, then $G = \G_1 = G'$.

{\it Inductive Step.} If $\G_2$ is a tree on $n$ vertices, 
then $\G$ is obtained from $\G_1$ by adding two subtrees 
of $\G_2$ with a total of $n-2$ edges. By Proposition \ref{tree}, 
$\Hr(\G) = \Hr(\G_1)\otimes \R(q^{n-2})$ 
regardless of the orientation of $e$.

If $\G_2$ is not a tree, let $e' \neq e \in \G_2$ be part of 
some cycle in $G_2$. Then 
$\Hr^i(G) = \Hr^{i-1}(G / e') \oplus \Hr^i(G - e')$, 
where $G / e'$ and $G - e'$ are obtained by gluing $\G_2 / e'$ and 
$\G_2 - e'$ respectively to $\G_1$ along $e$. Similarly, 
$\Hr^i(G') = \Hr^{i-1}(G' / e') \oplus \Hr^i(G' - e')$. By the inductive assumption, 
$\Hr^{i-1}(G / e') \cong \Hr^{i-1}(G' / e')$ and $\Hr^i(G - e') \cong \Hr^i(G' - e')$, 
and hence $\Hr^i(G) \cong \Hr^i(G')$.
\end{proof}

\begin{proposition}
The reduced cohomology sequence is an invariant of the matroid 
type of the graph.
\end{proposition}

\begin{proof}
The proof is by induction on the number of edges. There is only 
one graph with no edges, so the base of induction is vacuously true.

{\it Inductive Step.} If $\G$ is 2-connected, we are done by 
Proposition \ref{twist}. If $\G$ is not 2-connected, the 
removal of some vertex $v \in \G$ breaks $\G$ up into 
two connected components, $G_1$ and $G_2$. Adding $v$ back into 
$G_1$ and $G_2$ we get subgraphs $\G_1$ and $\G_2$ of 
$\G$, where $\G = \G_1 * \G_2$ (the vertex being $v$). 
Then the matroid type of $\G$ is the same as of the disjoint 
union of $\G_1$ and $\G_2$. On the other hand, 
$\Hr(\G) = \Hr(\G_1) \otimes \Hr(\G_2)$, 
which is also equal to the cohomology sequence of the disjoint 
union of $\G_1$ and $\G_2$.
\end{proof}
%-----------------------
%-----------------------
%5. Chromatic Polynomial
%-----------------------
%-----------------------

\section{Relationship to Standard Cohomologies and to the Chromatic Polynomial}\label{s:polynomial}

Here we derive the relationship between the reduced cohomologies of graphs 
and the main diagonal of the usual cohomologies. Then we describe the 
Poincar\'e polynomial for the reduced cohomologies in terms of the reduced 
chromatic polynomial.

\begin{proposition} For $C_n$ a cycle with $n$ vertices, 
$\Hr^i(C_n) = \R(q^{n-i-1})$ 
when $0 \leq i \leq n-2$; outside this range, 
$\Hr^n(C_n) = 0$.
\end{proposition}

\begin{proof}
We induct on $n$.

{\it Base Case.} If $n=1$, we have a loop, whose cohomologies 
are zero. If $n=2$, the graph has two vertices and two edges, 
so its cohomology sequence is $H^0(C_2) = \R(q)$ 
(and zero for the first and greater cohomology groups).

{\it Inductive Step.} If $i=0$, the zeroth cohomology group of any 
connected $n$-vertex graph is $\R(q^{n-1})$. For $1 \leq i$, 
$\Hr^i(C_n) \cong \Hr^{i-1}(C_n / e) \oplus \Hr^i(C_n - e)$. 
Now $C_n / e = C_{n-1}$, and $C_n - e = T_n$, a tree on $n$ 
vertices. By inductive assumption, 
$\Hr^{i-1}(C_{n-1}) = \R(q^{n-i-1})$ 
whenever $0 \leq i-1 \leq n-3$, i.e. $1 \leq i \leq n-2$; outside 
this range, $\Hr^{i-1}(C_{n-1}) = 0$. 
$\Hr^i(T_n) = \R(q^{n-1})$ if $i=0$, and $0$ 
otherwise. Adding the two, we get the statement of the proposition.
\end{proof}

\begin{definition}
If $\G$ is a graph, the Poincar\'e polynomial of $\G$ is a polynomial of two variables
$$\Rr_\G(t,q) := \sum_{i, j} t^iq^j\dim (\Hr^{i,j}(\G)).$$
\end{definition}

\begin{proposition}\label{deletion contraction}
Let $\G$ be a simple connected graph with $n$ vertices. Let $e$ be an edge of $\G$ 
that is not a bridge. Then the following 
\[
\Rr_\G(t,q) = t \Rr_{\G / e}(t,q) + \Rr_{\G - e}(t,q).
\]
\end{proposition}

\begin{proof} 
By Proposition \ref{direct sum}, $\dim(\Hr^{i,j}(\G)) = \dim(\Hr^{i-1,j}(\G/e)) + \dim(\Hr^{i,j}(\G-e))$.
Then
\[\begin{array}{lcl}
\Rr_\G(t,q) & = & \sum_{i, j} t^iq^j\dim (H^{i,j}(\G)) = \\
          & = &  \sum_{i, j} t^iq^j(\dim(\Hr^{i-1,j}(\G/e)) + \dim(\Hr^{i,j}(\G-e))) =\\
          & = &  t\sum_{i, j} t^{i-1}q^j\dim(\Hr^{i-1,j}(\G/e)) 
                                + \sum_{i, j} t^iq^j \dim(\Hr^{i,j}(\G-e)) = \\
          & = &  t \Rr_{\G / e}(t,q) + \Rr_{\G - e}(t,q).
\end{array}\]
\end{proof}

\begin{proposition}
Let $\G$ be a simple connected graph with $n$ vertices. Let $R^n_\G(t,q)$ be 
the homogeneous part of degree $n$ of the 
Poincar\'e polynomial for the usual cohomologies, and similarly for 
$R^{n-1}_\G(t,q)$. The Poincar\'e polynomial for the 
reduced cohomologies is
\[
 \Rr_\G (t,q) = 
\begin{cases} 
  \frac{1}{q} \left(R^n_\G (t,q)\left(1+\frac{t}{q}\right) - tq^{n-1}\right) = \vspace{10pt}\\ \hspace{2cm} 
    \frac{1}{q} R^n_\G(t,q) + R^{n-1}_\G(t,q) - q^{n-1},& \text{ if $\G$ is bipartite} \vspace{10pt}\\
  \frac{1}{q} R^n_\G(t,q)\left(1 + \frac{t}{q}\right)  = 
    \frac{1}{q} R^n_\G(t,q) + R^{n-1}_\G(t,q), & \quad\text{ otherwise}
\end{cases}
\]
\end{proposition}

\begin{proof}
The second set of equalities is a direct consequence of 
\cite[Theorem 5.2]{CCR}. We prove the first set, by induction 
on the number of edges.

First, we examine the cases of a tree and of a single 
odd-length-cycle graph. If $\G$ is a tree, 
$\Rr_\G(t,q) = q^{n-1}$. On the other hand, 
the non-reduced cohomologies of $\G$ are 
$H(\G) = \R(q^n) \oplus \R(q^{n-1})$ 
(see \cite[Example 28]{HGR}). Thus, $R^n_\G(t,q) = q^n$ and 
$R^{n-1}_\G(t,q) = q^{n-1}$. We observe
\[
q^{n-1} = \frac{1}{q}\left(q^n\left(1+\frac{t}{q}\right) - tq^{n-1}\right),
\]
as expected since trees are bipartite.

If $\G$ is a single cycle of length $n$ (odd), then 
\[
\Rr_\G(t,q) = q^{n-1} + q^{n-2}t + \dotsc + qt^{n-2}
\]
The $n$th degree homogeneous part of the 
non-reduced Poincar\'e polynomial, from \cite[Example 29]{HGR}, is
\[
R^n_\G(t,q) = q^n + q^{n-2}t^2 + \dotsc + q^3t^{n-3}
\]
The proposition follows by explicit computation.

Note that if $\G = \G_1 * \G_2$ where $\G_2$ is a tree 
on $n$ vertices, then both the reduced and the non-reduced cohomologies 
of $\G$ are computed by taking the respective cohomologies of $\G_1$ 
and multiplying by $\R(q^{n-1})$. Since adding a tree in this 
fashion preserves the bipartite or non-bipartite property, it also preserves 
the equality of polynomials above.

Now we proceed to the proper induction step. If $\G$ is not bipartite and 
contains more than one cycle (the one-cycle case was discussed above), then 
$\G$ contains some edge $e$ that is not a bridge and $\G/e$ is non-bipartite 
(Pick the smallest odd cycle 
$C$ of $\G$. We know that $\G$ has some other cycle, $C' \neq C$. 
Pick an edge $e \in C' - C$.) By construction, both $\G / e$ and $\G - e$ 
are both non-bipartite. Hence, by \cite[Theorem 5.5]{CCR} and by Proposition 
\ref{deletion contraction}, 
\[
R_\G(t,q) = t R_{\G / e}(t,q) + R_{\G - e}(t,q),
\]
\[
\Rr_\G(t,q) = t \Rr_{\G / e}(t,q) + \Rr_{\G - e}(t,q).
\]
By inductive assumption, the terms on the right-hand side 
satisfy the correct relations. Then
\[
\Rr_\G(t,q) = 
  \frac{t}{q} R^{n-1}_{\G/e}(t,q)\left(1 + \frac{t}{q}\right) + 
  \frac{1}{q} R^n_{\G - e}(t,q)\left(1 + \frac{t}{q}\right) = 
  \frac{1}{q} R^n_\G(t,q) \left(1 + \frac{t}{q}\right).
\]

If $\G$ is bipartite and contains a cycle, we take an edge $e$ contained 
in some (even-length) cycle. Then $\G - e$ will still be bipartite, but 
$\G / e$ will not be bipartite. Notice that in this case the same 
deletion-contraction relations for both standard and reduced homologies 
still hold. Therefore,
\begin{multline*}
\Rr_\G(t,q) = 
     \frac{t}{q} R^{n-1}_{\G/e}(t,q)\left(1 + \frac{t}{q}\right) + 
      \frac{1}{q} \left(R^n_{\G - e}(t,q)\left(1 + \frac{t}{q}\right) - tq^{n-1}\right)\\
= \frac{1}{q} \left(R^n_\G(t,q) \left(1 + \frac{t}{q}\right) - tq^{n-1}\right).
\end{multline*}
\end{proof}

In the standard case, \cite[Theorem 5.2]{CCR} derives $R^{n}_\G(-1,q)$ from 
$p_\G(q)$. The corresponding result for the reduced 
cohomologies is much simpler. The chromatic polynomial is the specification 
of the Poincar\'e polynomial at $t=-1$. The Poincar\'e polynomial for 
the reduced cohomologies is homogeneous of degree $n-1$ by Proposition 
\ref{diagonal}. Thus, $\Rr_\G(t,q)$ is completely 
determined by $\Rr_\G(-1,q)$. Specifically, 
$\Rr_\G(t,q) = (-t)^{n-1} p_\G\left(\frac{t-q}{t}\right)$.

\begin{remark}
One can prove that the reduced cohomologies over $\Z$ do not have any torsion. Hence 
all the results proven in this paper hold for cohomologies over $\Z$.
\end{remark}

\bigskip

\parbox{2in}
{\textbf{Michael Chmutov}

Department of Mathematics,

The Ohio State University, 

231 W. 18th Avenue, 

Columbus, Ohio 43210

\texttt{chmutov@mps.ohio-state.edu}}
\hspace{35pt}
\parbox{2in}
{\textbf{Elena Udovina}

Department of Mathematics,

Harvard University, 

One Oxford Street,

Cambridge, MA 02138

\texttt{eudovina@fas.harvard.edu}}

\end{document}